\documentclass[review]{elsarticle}

\usepackage{bm,amssymb,amsmath,mathrsfs}
\usepackage{amsthm}
\usepackage{lineno}
\usepackage[usenames]{color}

\usepackage{dcolumn}

\newcommand{\bq}{\begin{equation}}
\newcommand{\eq}{\end{equation}}
\newcommand{\bc}{\begin{center}}
\newcommand{\ec}{\end{center}}
\newcommand{\bit}{\begin{itemize}}
\newcommand{\eit}{\end{itemize}}
\newcommand{\ben}{\begin{enumerate}}
\newcommand{\een}{\end{enumerate}}

\theoremstyle{plain}

\newtheorem*{theorem*}{Theorem}

\begin{document}

\journal{(internal report CC23-12)}

\begin{frontmatter}

\title{Moments of the Poisson distribution of order $k$}

\author[cc]{S.~R.~Mane}
\ead{srmane001@gmail.com}
\address[cc]{Convergent Computing Inc., P.~O.~Box 561, Shoreham, NY 11786, USA}

\begin{abstract}
The factorial moments of the standard Poisson distribution are well known and are simple,
but the raw moments are considered to be more complicated (Touchard polynomials).
The present note presents a recurrence relation and an explicit combinatorial sum for the raw moments of the Poisson distribution of order $k$.
Unlike the standard Poisson distribution (the case $k=1$), for $k>1$ the structure of the raw and factorial moments have many similarities
and the raw moments are not more complicated (formally, at least) than the factorial moments.
We remark briefly on the central moments (i.e.~moments centered on the mean) of the Poisson distribution of order $k$.
\end{abstract}

\vskip 0.25in


\end{frontmatter}

\newpage
\setcounter{equation}{0}
It is well known that for a Poisson distribution with rate parameter $\lambda$, the mean and variance both equal $\lambda$
and the $n^{th}$ factorial moment is $M_{(n)}=\lambda^n$ for all $n\ge1$.
However, the expression for the raw moment $M_n$ is more complicated, and is given by a Touchard polynomial,
viz.~$M_0=1$ and for $n\ge1$ it is 
\bq
\label{eq:Poisson_rawmom}
M_n = \sum_{j=1}^k {n \brace j} \lambda^j \,.
\eq
Here ${n \brace j}$ is the Stirling number of the second kind.
The first few raw moments are
\begin{subequations}
\begin{align}
M_1 &= \lambda \,,
\\
M_2 &= \lambda^2 +\lambda \,,
\\
M_3 &= \lambda^3 +3\lambda^2 +\lambda \,,
\\
M_4 &= \lambda^4 +6\lambda^3 +7\lambda^2 +\lambda \,,
\\
M_5 &= \lambda^6 +10\lambda^4 +25\lambda^3 +15\lambda^2 +\lambda \,,
\\
M_6 &= \lambda^6 +15\lambda^5 +65\lambda^4 +90\lambda^3 +31\lambda^2 +\lambda \,.
\end{align}
\end{subequations}
The Poisson distribution of order $k$ \cite{PhilippouGeorghiouPhilippou} is a variant or extension of the standard Poisson distribution.
For $k=1$ it is the standard Poisson distribution.
It is an example of a compound Poisson distribution \cite{Adelson1966}.
Its probability mass function is \cite{PhilippouGeorghiouPhilippou} 
\bq
\label{eq:pmf}
P_n(k,\lambda) = e^{-k\lambda}\sum_{n_1+2n_2+\dots+kn_k=n} \prod_{j=1}^k \frac{\lambda^{n_j}}{n_j!} \,.
\eq
For $k=1$ this simplifies to the known result $P_n(\lambda) = e^{-\lambda}\lambda^n/n!$.
Less is known about its moments.
Let $M_n(k,\lambda)$ and $M_{(n)}(k,\lambda)$ respectively denote the $n^{th}$ raw and factorial moment of the Poisson distribution of order $k$.
Charalambides \cite{Charalambides1986} published an elegant analysis of statistical distributions of order $k$, including the Poisson distribution of order $k$.
In a recent note, \cite{Mane_Poisson_k_CC23_11}, the author presented the following combinatorial sum formula for the $n^{th}$ factorial moment.
It agrees with the formalism presented in \cite{Charalambides1986}.
\bq
\label{eq:facmom_orderk}
M_{(n)}(k,\lambda) = n!\sum_{n_1+2n_2+\dots+kn_k=n} \prod_{j=1}^k \frac{\lambda^{n_j}}{n_j!}\biggl[\binom{k+1}{j+1}\biggr]^{n_j} \,.
\eq

In this note, we present a combinatorial sum for the raw moment $M_n$ of the Poisson distribution of order $k$.
The formula has a similar structure to that for the factorial moments.
It is simplest to see this by employing the moment generating function and factorial moment generating function.
The probability generating function $G(z)$ (which is also the factorial moment generating function) and the moment generating function $M(t) = G(e^t)$
of the Poisson distribution of order $k$ are (\cite{Adelson1966}, see also \cite{Charalambides1986})
\begin{subequations}
\begin{align}
G(z,k,\lambda) &= e^{-k\lambda}e^{\lambda(z+z^2+\dots+z^k)} \,,
\\
M(t,k,\lambda) &= e^{-k\lambda}e^{\lambda(e^t+e^{2t}+\dots+e^{kt})} \,.
\end{align}
\end{subequations}
To simplify the derivations below, we introduce some useful polynomials.
Define the falling factorial sum $F_j(k,z)$ and power sum $S_j(k,t)$ as follows, for $j=1,2,\dots$
\begin{subequations}
\begin{align}
F_j(k,z) &= \sum_{s=j}^k s_{(j)}z^{s-j} = \sum_{s=j}^k s(s-1)\dots(s-j+1)z^{s-j} \,,
\\
S_j(k,t) &= \sum_{s=1}^k s^j e^{st} \,.
\end{align}
\end{subequations}
Also $F_0(k,z) = \sum_{s=1}^k z^s$ and $S_0(k,t) = \sum_{s=1}^ke^{st}$.
Note that $F_j(k,z)=0$ if $j>k$ but $S_j(k,t)\ne0$ for any value of $j$.
Also denote their values at $t=0$ and $z=1$ respectively as follows
\begin{subequations}
\begin{align}
\mathcal{F}_j(k) &= \sum_{s=j}^k s(s-1)\dots(s-j+1) &=&\; j!\binom{k+1}{j+1} \,,
\\
\mathcal{S}_j(k) &= \sum_{s=1}^k s^j
&=&\; \sum_{s=1}^k {j \brace s} \frac{(k+1)_{s+1}}{s+1}
&=&\; \sum_{s=1}^k {j \brace s} s!\binom{k+1}{s+1} \,.
\end{align}
\end{subequations}
Then observe that
\bq
G(z,k,\lambda) = e^{-k\lambda}e^{\lambda F_0(k,z)} \,,
\qquad\qquad
M(t,k,\lambda) = e^{-k\lambda}e^{\lambda S_0(k,t)} \,.
\eq
Note the derivatives
\bq
\frac{dF_j(k,z)}{dz} = F_{j+1}(k,z) \,,
\qquad\qquad
\frac{dS_j(k,t)}{dt} = S_{j+1}(k,t) \,.
\eq
We omit explicit mention of the function arguments below.
The first few derivatives are given as follows
\begin{subequations}
\begin{align}
\frac{dG}{dz} &= e^{-k\lambda}e^{\lambda F_0} F_1\lambda  \,,
&
\frac{dM}{dt} &= e^{-k\lambda}e^{\lambda S_0} S_1\lambda  \,,
\\
\frac{d^2G}{dz^2} &= e^{-k\lambda}e^{\lambda F_0} (F_1^2\lambda ^2 + F_2\lambda ) \,,
&
\frac{d^2M}{dt^2} &= e^{-k\lambda}e^{\lambda S_0} (S_1^2\lambda^2 + S_2\lambda ) \,,
\\
\frac{d^3G}{dz^3} &= e^{-k\lambda}e^{\lambda F_0} (F_1^3\lambda^3 +3F_1F_2\lambda^2 +F_3\lambda) \,,
&
\frac{d^3M}{dt^3} &= e^{-k\lambda}e^{\lambda F_0} (S_1^3\lambda^3 +3S_1S_2\lambda^2 +S_3\lambda) \,,
\\
\frac{d^4G}{dz^4} &= e^{-k\lambda}e^{\lambda F_0} (F_1^4\lambda^4 +6F_1^2F_2\lambda^3 
&
\frac{d^4M}{dt^4} &= e^{-k\lambda}e^{\lambda S_0} (S_1^4\lambda^4 +6S_1^2S_2\lambda^3
\nonumber\\
& \qquad\qquad\qquad
+(4F_1F_3 +3F_2^2)\lambda^2 +F_4\lambda) \,,
&& \qquad\qquad\qquad
+(4S_1S_3 +3S_2^2)\lambda^2 +S_4\lambda) \,.
\end{align}
\end{subequations}
Set $z=0$ to obtain the factorial moment and $t=1$ for the raw moment.
We obtain 
\begin{subequations}
\begin{align}
M_{(1)} &= \mathcal{F}_1\lambda \,,
& M_1 &= \mathcal{S}_1\lambda \,,
\\
M_{(2)} &= \mathcal{F}_1^2\lambda^2 + \mathcal{F}_2\lambda \,,
&
M_2 &= \mathcal{S}_1^2\lambda^2 + \mathcal{S}_2\lambda \,,
\\
M_{(3)} &= \mathcal{F}_1^3\lambda^3 +3\mathcal{F}_1\mathcal{F}_2\lambda^2 +\mathcal{F}_3\lambda \,,
&
M_3 &= \mathcal{S}_1^3\lambda^3 +3\mathcal{S}_1\mathcal{S}_2\lambda^2 +\mathcal{S}_3\lambda \,,
\\
M_{(4)} &= \mathcal{F}_1^4\lambda^4 +6\mathcal{F}_1^2\mathcal{F}_2\lambda^3 +(4\mathcal{F}_1\mathcal{F}_3 +3\mathcal{F}_2^2)\lambda^2 +\mathcal{F}_4\lambda \,,
&
M_4 &= \mathcal{S}_1^4\lambda^4 +6\mathcal{S}_1^2\mathcal{S}_2\lambda^3 +(4\mathcal{S}_1\mathcal{S}_3 +3\mathcal{S}_2^2)\lambda^2 +\mathcal{S}_4\lambda \,.
\end{align}
\end{subequations}
The two sets of expressions are the same, just replace $\mathcal{F}_j$ by $\mathcal{S}_j$.
It is easily seen that this substitution applies for all the higher derivatives, i.e.~the higher moments.
The recurrence relations for the probability mass function $P_n$, factorial moment $M_{(n)}$ and raw moment $M_n$ are respectively
\begin{subequations}
\begin{align}
\label{eq:rec_Pn_me}    
P_n &= \frac{\lambda}{n} \sum_{j=1}^{\min\{n,k\}} jP_{n-j} \,,
\\
\label{eq:rec_facmom_me}    
M_{(n)} &= \lambda \sum_{j=1}^{\min\{n,k\}} \binom{n-1}{j-1}\mathcal{F}_jM_{(n-j)} \,,
\\
\label{eq:rec_rawmom_me}    
M_n &= \lambda \sum_{j=1}^n \binom{n-1}{j-1}\mathcal{S}_jM_{n-j} \,.
\end{align}
\end{subequations}
Note the following:
\begin{enumerate}
\item
  First, eq.~\eqref{eq:rec_Pn_me} was derived in \cite{Adelson1966} for compound Poisson distributions in general,
  and for the Poisson distribution of order $k$ in \cite{GeorghiouPhilippouSaghafi} and \cite{KostadinovaMinkova2013}.
\item
  Next, eq.~\eqref{eq:rec_facmom_me} was derived by Charalambides (eq.~(4.7) in \cite{Charalambides1986}, with differences of notation)
  and eq.~\eqref{eq:rec_rawmom_me} is obtained by analogy with eq.~\eqref{eq:rec_facmom_me},  
  replacing $\mathcal{F}_j$ by $\mathcal{S}_j$.
\item
  The upper limit for the sum in eq.~\eqref{eq:rec_Pn_me} cuts off at $\min\{n,k\}$ because $P_j=0$ for $j<0$.
\item
  The upper limit for the sum in eq.~\eqref{eq:rec_facmom_me} cuts off at $\min\{n,k\}$ because $\mathcal{F}_j=0$ for $j>k$.
\item
  The upper limit for the sum in eq.~\eqref{eq:rec_rawmom_me} is $n$ because $\mathcal{S}_j\ne0$ for all for $j\ge0$.
\end{enumerate}
The combinatorial sum expressions for $P_n$, $M_{(n)}$ and $M_n$ are respectively
\begin{subequations}
\begin{align}
\label{eq:Pn_sol_me}    
P_n &= e^{-k\lambda} \sum_{n_1+2n_2+\dots+kn_k=n} \prod_{j=1}^k \frac{\lambda^{n_j}}{n_j!} \,,
\\
\label{eq:facmom_sol_me}    
M_{(n)} &= n! \sum_{n_1+2n_2+\dots+kn_k=n} \prod_{j=1}^k \frac{\lambda^{n_j}}{n_j!} \biggl(\frac{\mathcal{F}_j}{j!}\biggr)^{n_j}
&= n! \sum_{n_1+2n_2+\dots+kn_k=n} \prod_{j=1}^k \frac{\lambda^{n_j}}{n_j!}\biggl[\binom{k+1}{j+1}\biggr]^{n_j} \,,
\\
\label{eq:rawmom_sol_me}    
M_n &= n! \sum_{n_1+2n_2+\dots+kn_k=n} \prod_{j=1}^k \frac{\lambda^{n_j}}{n_j!} \biggl(\frac{\mathcal{S}_j}{j!}\biggr)^{n_j} \,.
\end{align}
\end{subequations}
Note the following:
\begin{enumerate}
\item
First, eq.~\eqref{eq:Pn_sol_me} is the same as eq.~\eqref{eq:pmf}
and eq.~\eqref{eq:facmom_sol_me} is the same as eq.~\eqref{eq:facmom_orderk}.
They are reproduced here for ease of reference, and to display the similarity of structure with eq.~\eqref{eq:rawmom_sol_me}.
\item
For the standard Poisson distribution ($k=1$), the expression for $P_n$ contains only a single term $e^{-\lambda}\lambda^n/n!$ 
and the expression for $M_{(n)}$ contains only a single term $\lambda^n$, 
whereas the expression for $M_n$ is a polynomial with $k$ terms (see eq.~\eqref{eq:Poisson_rawmom}).
\item
However, for $k>1$, the expression for the raw moment $M_n$ is not more complicated (formally, at least) than for the factorial moment $M_{(n)}$.
Admittedly, the expression for $\mathcal{F}_j$ is simpler than that for $\mathcal{S}_j$.
\end{enumerate}

\newpage
We remark on the central moments (i.e.~moments centered on the mean), denoted by $\tilde{M}_n = \mathbb{E}[(X-\mathbb{E}[X])^n$.
The mean of the Poisson distribution of order $k$ is (\cite{Adelson1966}, see also \cite{PhilippouMeanVar})
\bq
\mu(k,\lambda) = (1+\dots+k)\lambda = \frac{k(k+1)}{2}\lambda = \mathcal{S}_1\lambda \,.
\eq
The value of $\tilde{M}_2$ is the variance (\cite{Adelson1966}, see also \cite{PhilippouMeanVar})
\bq
\sigma^2(k,\lambda) = (1^2+\dots+k^2)\lambda = \frac{k(k+1)(2k+1)}{6}\lambda = \mathcal{S}_2\lambda \,.
\eq
For the standard Poisson distribution ($k=1$), both the mean and the variance equal $\lambda$.
Then for all $n\ge2$, the formula for $\tilde{M}_n$ is simply that for the noncentral raw moment $M_n$ {\em without} the terms in $\mathcal{S}_1$ (see eq.~\eqref{eq:rawmom_sol_me})
\bq
\label{eq:centeredmom_sol_me}    
\tilde{M}_n = n! \sum_{2n_2+\dots+kn_k=n} \prod_{j=2}^k \frac{\lambda^{n_j}}{n_j!} \biggl(\frac{\mathcal{S}_j}{j!}\biggr)^{n_j} \,.
\eq
From eq.~\eqref{eq:rec_rawmom_me}, the recurrence relation for the centered moments is (with $\tilde{M}_0=1$ and $\tilde{M}_1=0$, by definition)
\bq
\label{eq:rec_centeredmom_me}    
\tilde{M}_n = \lambda \sum_{j=2}^n \binom{n-1}{j-1}\mathcal{S}_j\tilde{M}_{n-j} \,.
\eq
The first few central moments of the Poisson distribution of order $k$ are
\begin{align}
\tilde{M}_2 &= \mathcal{S}_2\lambda \,,
&
\tilde{M}_3 &= \mathcal{S}_3\lambda \,,
\nonumber\\
\tilde{M}_4 &= 3\mathcal{S}_2^2\lambda^2 +\mathcal{S}_4\lambda \,,
&
\tilde{M}_5 &= 10\mathcal{S}_2\mathcal{S}_3\lambda^2 +\mathcal{S}_5\lambda \,,
\nonumber\\
\tilde{M}_6 &= 15\mathcal{S}_2^3\lambda^3 +(15\mathcal{S}_2\mathcal{S}_4+10\mathcal{S}_3^2)\lambda^2 +\mathcal{S}_6\lambda \,,
&
\tilde{M}_7 &= 105\mathcal{S}_2^2\mathcal{S}_3\lambda^3 +(21\mathcal{S}_2\mathcal{S}_5+35\mathcal{S}_3\mathcal{S}_4)\lambda^2 +\mathcal{S}_7\lambda \,.
\end{align}
Note $\tilde{M}_n$ is a polynomial in $\lambda$ of degree $\lfloor(n/2)\rfloor$
with no constant term and all powers from $\lambda$ up to $\lambda^{\lfloor(n/2)\rfloor}$.
The coefficient of $\lambda$ is $\mathcal{S}_n$.
For the standard Poisson distribution ($k=1$), then $\mathcal{S}_j=1$ for all $j$ and the above expressions simplify to the known results.
For $k=1$, we write $\mu_n$ in place of $\tilde{M}_n$:
\begin{align}
\mu_2 &= \lambda \,,
&
\mu_3 &= \lambda \,,
\nonumber\\
\mu_4 &= 3\lambda^2 +\lambda \,,
&
\mu_5 &= 10\lambda^2 +\lambda \,,
\nonumber\\
\mu_6 &= 15\lambda^3 +25\lambda^2 +\lambda \,,
&
\mu_7 &= 105\lambda^3 +56\lambda^2 +\lambda \,.
\end{align}

In conclusion, the papers by Adelson \cite{Adelson1966} and Charalambides \cite{Charalambides1986}
lay an elegant foundation for the factorial moments of the Poisson distribution of order $k$.
This was extended to treat the raw and centered moments.
This note displays explicit combinatorial sums for the above quantities.
(The expression for the probability mass function was published in \cite{PhilippouGeorghiouPhilippou}
and the expression for the factorial moment was published in \cite{Mane_Poisson_k_CC23_11}, in agreement with formulas in \cite{Charalambides1986}.)
There is a consistent pattern in the combinatorial sums for all four of the above quantities, also for their recurrence relations.


\end{document}